\documentclass[11pt]{amsart}
\pdfoutput=1
\usepackage{amssymb}
\usepackage{amsmath}
\usepackage{amsfonts}
\usepackage[usenames]{color}
\usepackage{graphicx}
\usepackage{array}
\usepackage{psfrag}
\usepackage{color}

\makeatletter
 
 \@addtoreset{equation}{section}
\makeatother

\newtheorem{theorem}{Theorem}[section]

\theoremstyle{definition}

\theoremstyle{remark}

\numberwithin{equation}{section}

\newcommand{\R}{\mathbb{R}}

\begin{document}

\title{Braiding Knots in Contact 3-manifolds}

\author{Elena Pavelescu}
\address{ Rice University,  Department of Mathematics - MS 136, Houston, TX 77005}
\email{Elena.Pavelescu@rice.edu}



\begin{abstract}
We show that a transverse link in a contact structure
supported by an open book decomposition can be transversely
braided. We also generalize Markov's theorem on when the closures of
two braids represent (transversely) isotopic links.
\end{abstract}
\maketitle
\section{Introduction}

In ~\cite{Al}, Alexander proved that any knot in $\mathbb{R}^3$ can
be braided about the $z$-axis. In \cite{Be}, a paper which marked
the start of modern contact topology and led to the Bennequin
inequality and eventually to Eliashberg's definition of tightness,
Bennequin proved the transverse case for $(\mathbb{R}^3,\xi_{std})$.
In this paper, we generalize Bennequin's result to any closed
oriented 3-dimensional manifold $M$, by looking at an open book
decomposition for $M$ together with its supported contact structure.
After completing the proof of this theorem, it was brought to our
attention that the same result was independently obtained through different means by
Mitsumatsu and Mori in \cite{Mi}.

Bennequin used braid theory in $\mathbb{R}^3$ to show the standard
contact structure was tight. Since then, powerful analytic machinery
like holomorphic curves or Seiberg-Witten theory has been used to prove tightness. By studying braids
in general open books it is hoped that tight contact structures can be better
understood from a purely topological/combinatorial perspective. The
results presented in this paper constitute a first step in this
program. While these results are purely geometric, it would be nice
to have an algebraic description of a generalized $''$braid group$''$.
In particular, one needs to define the stabilization moves with
respect to different binding components and account for the
fundamental group of the page. We will explore this algebraic
perspective in a subsequent paper.

In \cite{Ma}, Markov gave an equivalent condition for the closures
or two braids in $\mathbb{R}^3$ to be isotopic as links. This is the
case if and only if the two braids differ by conjugations in the
braid group and positive and negative Markov moves. In \cite{OrSh},
Orevkov and Shevchishin proved the transversal case for
$(\mathbb{R}^3,\xi_{std})$. A different proof was independently
obtained by Wrinkle in \cite{Wr}. We generalize Markov's theorem to
any closed oriented 3-dimensional manifold. We prove the transverse
case and recover the topological case previously proved in \cite{Sk}
and \cite{Su}.

In \cite{KP}, the author together with Keiko Kawamuro used the braiding results presented in this paper to find combinatorial self-liking formulae for null-homologous braids in annulus open book decompositions and pants open book decompositions.  These formulae extend Bennequin's self-linking formula for a braid in the standard contact 3-sphere.

\section*{Acknowledgements}
The author thanks John Etnyre for all his support and guidance. 
She would also like to thank Evgeny Volkov for taking an interest in these results and pointing out the sign mistake in Theorem 3.2.
She thanks Keiko Kawamuro for all the helpful conversations.

\section{Contact Structures and Open Book Decompositions}

Let M be a compact, oriented 3-manifold, and $\xi$ an oriented 2-plane field on $M$.  We say $
\xi$ is a \emph{contact structure} on $M$ if $\xi =\ker \alpha$ for some global
1- form $\alpha$ satisfying $\alpha\wedge d\alpha\ne 0.$  In this paper we assume
that the manifold $M$ is oriented and that the contact structure $\xi$ is oriented
and positive (i.e. $\alpha\wedge d\alpha> 0$).

A \textit{contactomorphism} between two contact manifolds $(M_1,
\xi_1)$ and $(M_2, \xi_2)$ is a diffeomorphism $\phi :M_1\rightarrow
M_2$ such that $\phi_*\xi_1=\xi_2$.

On $\mathbb{R}^3$, consider the two contact structures $\xi_1$ and
$\xi_2$ given by the 1-forms $\alpha_1= dz-ydx$ and
$\alpha_2=dz+r^2d\theta$ (given in cylindrical coordinates). Then
$(\R^3,\xi_1)$ and $(\R^3,\xi_2)$  are contactomorphic and we are
going to refer to both $\xi_1$ and $\xi_2$ as the standard contact
structure on $\R^3$, $\xi_{std}$.

\begin{theorem} ( Gray )\label{Gr}
Let $\{\xi_t\}_{t\in [0,1]}$ be a family of contact structures on a
manifold $M$ that differ on a compact set $C\subset$ int$(M)$. Then
there exists an isotopy $\psi_t:M\rightarrow M$ such that
\begin{enumerate}
\item[i)] $(\psi_t)_*\xi_1$=$\xi_t$ \item[ii)] $\psi_t$ is the
identity outside of an open neighborhood of $C$.
\end{enumerate}
\label{thm:gray}
\end{theorem}

While the proof of this result is well known, we sketch it here as
we will need elements of it in later arguments. For more details see
\cite{Ge} (p. 61).
 
\begin{proof}
We are going to look for $\psi_t$ as the flow of
a vector field $X_t$. If $\xi_t=\ker \alpha_t$, then $\psi_t$ has to
satisfy
\begin{equation}\psi_t^*\alpha_t=\lambda_t\alpha_0,
\end{equation} for some non-vanishing function
$\lambda_t:M\rightarrow \mathbb{R}$. By taking the derivative with
respect to $t$ on both sides and rearranging the terms we get
\begin{equation}\label{eq1}
\psi_t^*(\frac{d\alpha_t}{dt}+\mathcal{L}_{X_t}\alpha_t)=\frac{d\lambda_t}{dt}\alpha_0
= \frac{d\lambda_t}{dt} \frac{1}{\lambda_t}\psi^*\alpha_t.
\end{equation}
This is equivalent to

\begin{equation}\label{eq2}
\psi_t^*(\frac{d\alpha_t}{dt}+d(\iota_{X_t}\alpha_t)+
\iota_{X_t}d\alpha_t)=\psi_t^*(h_t\alpha_t)
\end{equation}

for $$h_t=\frac{d}{dt}(\log \lambda_t)\circ \psi_t^{-1}.$$

If $X_t$ is chosen in $\xi_t$ then $\iota_{X_t}\alpha_t=0$ and
(\ref{eq2}) becomes
\begin{equation}\label{eq3}\frac{d\alpha_t}{dt}+\iota_{X_t}d\alpha_t
=h_t\alpha_t\end{equation}

Applying (\ref{eq3}) to the Reeb vector field of $\alpha_t$,
$v_{\alpha_t}$ (that is, the unique vector field $v_t$ such that
$\alpha_t(v_t)=1$ and $d\alpha_t(v_t, \cdot)=0$), we find $h_t =
\frac{d\alpha_t}{dt}(v_{\alpha_t})$ and $X_t$ given by
\begin{equation}\label{eq4}\iota_{X_t}d\alpha_t = h_t\alpha_t
-\frac{d\alpha_t}{dt}\end{equation}

The form $d\alpha_t$ gives an isomorphism
\begin{center}
$\Gamma(\xi_t) \rightarrow \Omega^1_{\alpha_t}$
\\ $v\mapsto \iota_vd\alpha_t$
\end{center}
where $\Gamma(\xi_t)= \{v|v\in \xi_t\}$ and $\Omega_{\alpha_t}^1
=\{$1-forms $\beta|\beta(v_t)=0\}$, and thus $X_t$ is uniquely
determined by (\ref{eq3}). By construction, the flow of $X_t$ is the desired $\psi_t$.
For the subset of $M$ where the $\xi_t$'s agree we choose the
$\alpha_t$'s to agree. This implies $\frac{d\alpha_t}{dt}=0$,
$h_t=0$ and $X_t=0$ and all equalities hold.
\end{proof}

In a contact manifold $(M, \xi)$, an oriented arc $\gamma \subset M$
is called \textit{transverse} if for all $p\in\gamma$ and $\xi_p$
the contact plane at $p$, $T_p\gamma \pitchfork \xi_p$ and
$T_p\gamma$ intersects $\xi_p$ positively.  Moreover, if $\gamma$ is a closed curve then it is
called a \textit{transverse knot}.

An \textit{open book decomposition} of $M$ is a pair (L, $\pi$)
where\begin{enumerate}
\item[(1)] $L$ is an oriented link in $M$ called \textit{the binding}
of the open book
\item[(2)] $\pi : M\smallsetminus L\to S^1$ is a
fibration whose fiber, $\pi^{-1}(\theta)$, is the interior of a
compact surface $\Sigma\subset M$ such that $\partial \Sigma = L$,
for all $ \theta\in S^1$. The surface $\Sigma$ is called \textit{the
page} of the open book.
\end{enumerate}

Alternatively, an open book decomposition of a 3-manifold $M$
consists of a surface $\Sigma$, with boundary, together with a
diffeomorphism $\phi: \Sigma\rightarrow \Sigma$, with $\phi$ the
identity near $\partial \Sigma$, such that
$$M=(\Sigma\times[0,1]/\sim) \cup_f \coprod_i S^1\times D^2$$ where
$$(x,1)\sim(\phi(x),0).$$ Note that $$\partial(\Sigma\times [0,1]/\sim) = \coprod_i T_i^2,$$
with each torus $T_i^2$ having a product structure $S^1\times
[0,1]/\sim$. Let $\lambda_i =\{pt\}\times [0,1]/\sim$, $\lambda_i\in
T_i^2$. The gluing diffeomorphism used to construct $M$ is defined
by
\begin{center}
$$f:\partial(\coprod_i S^1\times D^2)\rightarrow \partial
(\coprod_i T_i^2)$$\\
$\{pt\}\times \partial D^2 \rightarrow \lambda_i.$
\end{center} The map $\phi$ is called \textit{the monodromy} of
the open book. See \cite{Et} for more details.
\begin{theorem} (Alexander, \cite{Al1}) Every closed oriented 3-manifold has an open book
decomposition.\end{theorem}

A contact structure $\xi$ on $M$ is said to be \textit{supported} by
an open book decomposition $(\Sigma, \phi)$ of $M$ if $\xi$ can be
isotoped through contact structures so that there exists a 1-form
$\alpha$ for $\xi$ such that
\begin{enumerate}
\item[(1)] $d\alpha$ is a positive area form on each page
\item[(2)] $\alpha(v)>0$ for all $v\in TL$ that induce the
orientation on $L$.
\end{enumerate}

\begin{theorem}(Thurston, Winkelnkemper, \cite{ThWi}) Every open
book decomposition $(\Sigma, \phi)$ supports a contact structure
$\xi_{\phi}$. \label{Thm:ThWi}
\end{theorem}
We sketch this well-known proof as we need the details in later arguments.

\begin{proof}
Let $$M=(\Sigma\times[0,1]/\sim) \cup_f \coprod_i S^1\times D^2$$
given as before. We first construct a contact structure on
$\Sigma\times[0,1]/\sim$ and then we extend it in a neighborhood
of the binding.

In the neighborhood $N=S^1\times D^2$ of each component of the
binding consider coordinates $(\psi, x,\theta)$ such that ($\psi,
x)$ are coordinates on the page with $\psi$ the coordinate along the
binding and $d\theta$ and $\pi^*d\theta$ agree, where $\pi^*d\theta$
is the pullback through $\pi: M\setminus L \rightarrow S^1$ of the
coordinate on $S^1$. Let $\lambda$ be a 1-form on the page which is
an element of the set (it is easy to check that this set is not
empty)
\begin{eqnarray*}
S=\{1-\textrm{forms }\lambda \textrm{ such that} & : & d\lambda
\textrm{ is a volume form on } \Sigma{}\nonumber{}\nonumber
\\ & & {} \lambda =(1+x)d\psi \textrm{ near }
\partial\Sigma=L\}
\end{eqnarray*}
On $\Sigma \times [0,1]$ take $\tilde{\lambda}= (1-\theta)\lambda
+\theta(\phi^*\lambda)$ and consider the 1-form
$\alpha_K=\tilde{\lambda}+Kd\theta$. For sufficiently large $K$,
$\alpha_K$ is a contact form and it descends to a contact form on
$\Sigma\times[0,1]/\sim$. To extend this form on the solid tori
neighborhood of the binding we pull back $\alpha$ through the gluing
map $f$ and get  \begin{center} $\alpha_f=
Kd\theta-(x+\epsilon)d\psi.$
\end{center}

\begin{figure}[htpb!]
\begin{center}
\begin{picture}(260, 144)
\put(0,0){\includegraphics{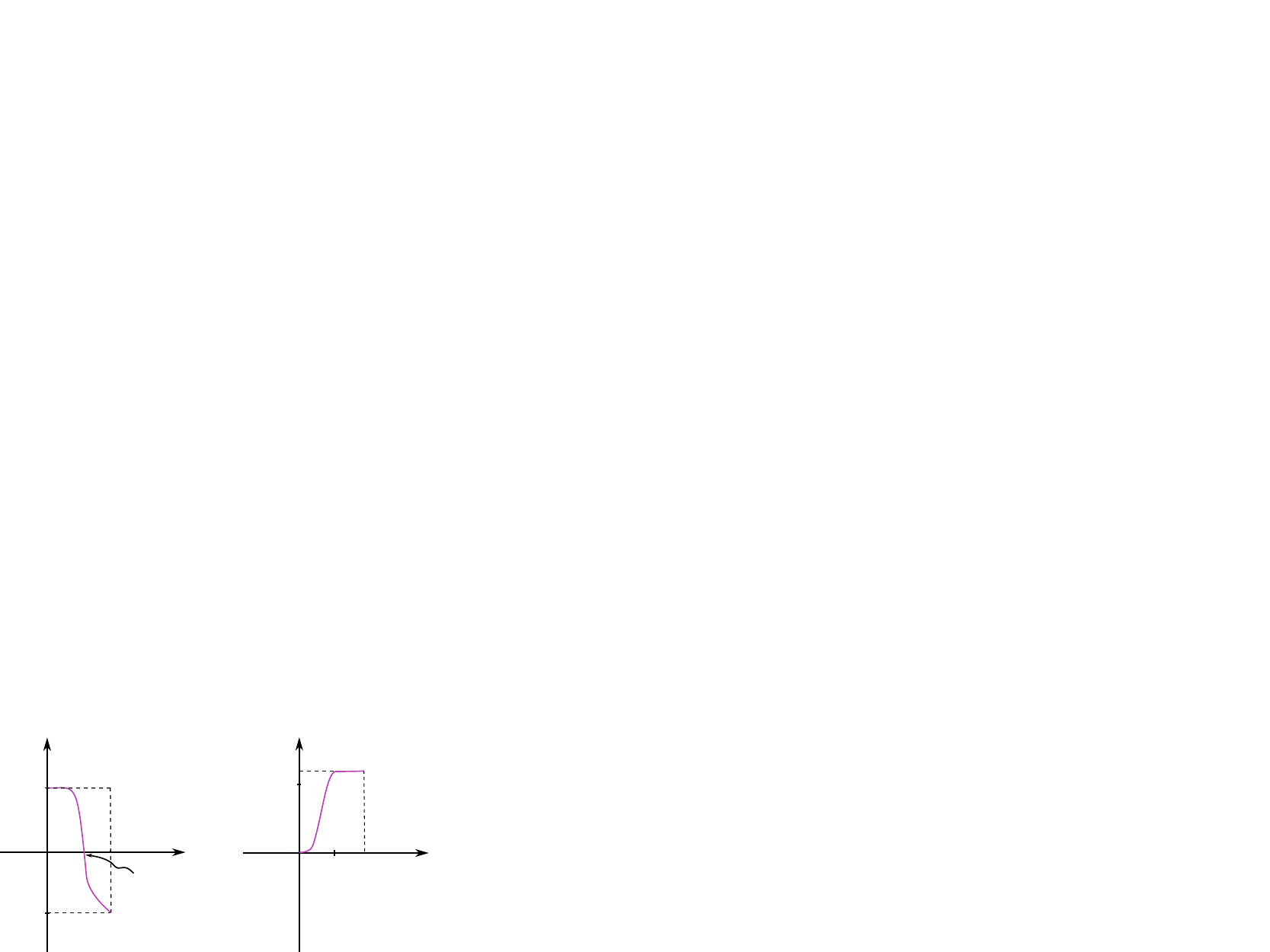}}
\put(-5, 22){$-1-\epsilon$}
\put(55, 80){$h$}
\put(203, 80){$g$}
\put(265,53){$x$}
\put(110,53){$x$}
\put(176,111){$K$}
\put(22,100){$1$}
\put(89,47){$\delta_{h}$}
\put(208,52){$\delta_{g}$}
\end{picture}
\caption{$h$ and $g$ functions}\label{hg-functions}
\end{center}
\end{figure}
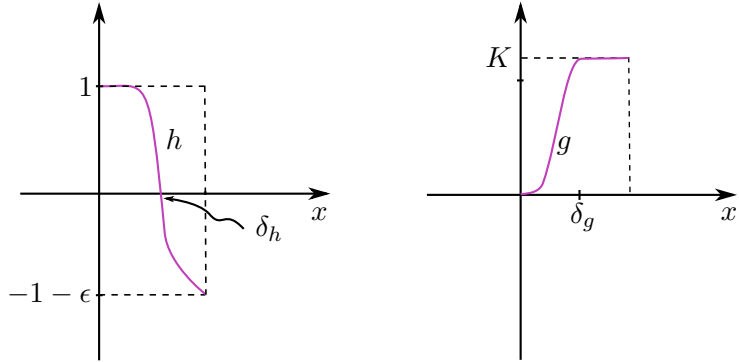

We are looking to extend this form on the entire $S^1\times D^2$ to
a contact form of the form $h(x)d\psi+g(x)d\theta$. This is possible
if there exist functions $h,g:[0,1]\rightarrow \mathbb{R}^3$ such
that:
\begin{enumerate}\item[(1)]$h(x)g'(x)-h'(x)g(x)>0$ (given by the
contact condition)\item[(2)]$h(x)=1$ near $x=0$,
$h(x)=-(x+\epsilon)$ near $x=1$ \item[(3)]$g(x)=r^2$ near $x=0$,
$g(x)=K$ near $x=1$.
\end{enumerate}

The two functions $h$ and $g$ described in Figure ~\ref{hg-functions}
work for our purpose. 
The conditions (2) and (3) are obviously
satisfied, and if $\delta_h$ and $\delta_g$ are such that $h(x)<0$ on
$[\delta_h, 1]$ and $g(x)=K$ on $[\delta_g,1]$, then (1) is satisfied
as long as $\delta_g<\delta_h$.

\end{proof}

The following correspondence, which is one of the central results in
contact geometry, relates open book decompositions and contact
structures.
\begin{theorem}(Giroux, \cite{Gi})
Let M be a closed, oriented 3-manifold. Then there is a one to one
correspondence

 $$\left\{
\begin{array}{l}
\mbox{Contact structures $\xi$ on $M^3$} \\
\mbox{up to contact isotopy}
\end{array}
\right\} \stackrel{1-1}{\longleftrightarrow} \left\{
\begin{array}{l}
\mbox{Open book decompositions $(\Sigma, \phi)$} \\
\mbox{of $M^3$ up to positive stabilization}
\end{array}
\right\}.
$$
\end{theorem}
See \cite{Et} for details.


\section{Braiding Knots in contact 3-manifolds}
In this section we generalize the following theorem:
\begin{theorem} (Bennequin, \cite{Be})
Any transverse link  $K$ in $(\mathbb{R}^3, \xi_{std})$ is
transversely isotopic to a link braided about the $z$-axis.
\label{thm:bennequin}
\end{theorem}

Let $(L, \pi)$ be an open book decomposition for $M$. A link
$K\subset M$ is said to be \textit{braided about} $L$ if $K$ is
disjoint from $L$ and there exists a parametrization of $K$, $f:
\coprod S^1 \rightarrow M$, such that if $\theta$ is the coordinate
on each $S^1$ then $\frac{d}{d\theta}(\pi\circ f)>0$, $\forall
\theta$. We call \emph{bad arcs} of $K$ those arcs where this
condition is not satisfied.

\begin{theorem} Suppose $(L, \pi)$ is an open book
decomposition for the 3-manifold $M$ and $\xi$ is supported by $(L,
\pi)$. Let $K$ be a transverse link in $M$. Then $K$ can be
transversely isotoped to a braid.\label{Alexander}\end{theorem}

\begin{proof}
 The idea of the proof is to find a family of diffeomorphisms
of $M$ keeping each page of the open book setwise fixed and sending
the parts of the link where the link is not braided into a
neighborhood of the binding. A neighborhood of the binding is
contactomorphic to a neighborhood of the $z$-axis in
$(\mathbb{R}^3,\xi_{std})/_{z\sim z+1}$ and there the link can be
braided, according to Theorem \ref{thm:bennequin}. 
What we really use is that Theorem \ref{thm:bennequin} works not only for links, but also  
for arcs with good ends. 

In the neighborhood $N=S^1\times D^2$ of each component of the
binding consider coordinates $(\psi, x,\theta)$ and let $\lambda\in
S$ as described in Theorem \ref{Thm:ThWi}. On $\Sigma \times [0,1]$
take $\tilde{\lambda}= (1-\theta)\lambda +\theta(\phi^*\lambda)$ and
consider the family of 1-forms given by
$\alpha_t=\tilde{\lambda}-K\frac{1}{t}d\theta$, where  $t\in[-1,0)$
and $K$ is a large constant.

This family of 1-forms descends to a family of 1-forms on
$\Sigma\times [0,1]/\sim$ (because $\phi^*(\alpha_t|_{\Sigma\times
\{0\}})=\alpha_t|_{\Sigma\times\{1\}}$). Both $\xi_{-1}=ker(\alpha_{-1})$ and
$\xi$ are contact structures supported by $(L,\pi)$  and by Giroux's
correspondence they are isotopic. Therefore, without loss of
generality, we may assume $\xi =ker(\alpha_{-1})$.

For large enough $K$, the family of 1-forms $\{\alpha_t\}_t$ is a
family of contact 1-forms as:
$$\alpha_t\wedge d\alpha_t=
(\tilde{\lambda}-K\frac{1}{t}d\theta)\wedge (d\tilde{\lambda})=
\tilde{\lambda}\wedge d\tilde{\lambda} - K\frac{1}{t}d\theta \wedge
d\tilde{\lambda}>0$$ Note that $d\tilde{\lambda}$ is an area form on
the page while $d\theta$ vanishes on the page and is positive on the
positive normal to the page. This implies that the term which is subtracted
 is always negative ( $t\in [-1,0)$) and therefore $\alpha_t$ is a contact
form for sufficiently large $K$. We want to extend this family to
the whole $M$, so we need to patch in the solid tori neighborhood of
the binding. In the neighborhood of the binding we have $\alpha_t=(1+x)d\psi
- K\frac{1}{t}d\theta$. We pull back this family through the map $f$
used to glue the solid tori in the definition of the open book. With
our chosen coordinates we have
 $$f(\psi, x, \theta)=(x-1+\epsilon, -\psi, \theta)$$
 and the pullback family
 $$\alpha_{f,t}= -(x+\epsilon)d\psi- K\frac{1}{t}d\theta.$$

We are looking to extend these forms on the entire $S^1\times D^2$ to
a family of the form $h_t(x)d\psi+g_t(x)d\theta$. This is possible
if there exist functions $h_t,g_t:[0,1]\rightarrow \R$, for
$t\in[-1,0)$ such that:
\begin{enumerate}\item[i)]$h_t(x)g_t'(x)-h_t'(x)g_t(x)>0, t\in[-1,0)$ (given by the
contact condition)\item[ii)]$h_t(x)=1$ near $x=0$,
$h_t(x)=-(x+\epsilon)$ near $x=1$ \item[iii)]$g_t(x)=x^2$ near
$x=0$, $g_t(x)=-\frac{K}{t}$ near $x=1$.
\end{enumerate}

\begin{figure}[htpb]
\begin{center}
\begin{picture}(260, 144)
\put(0,0){\includegraphics{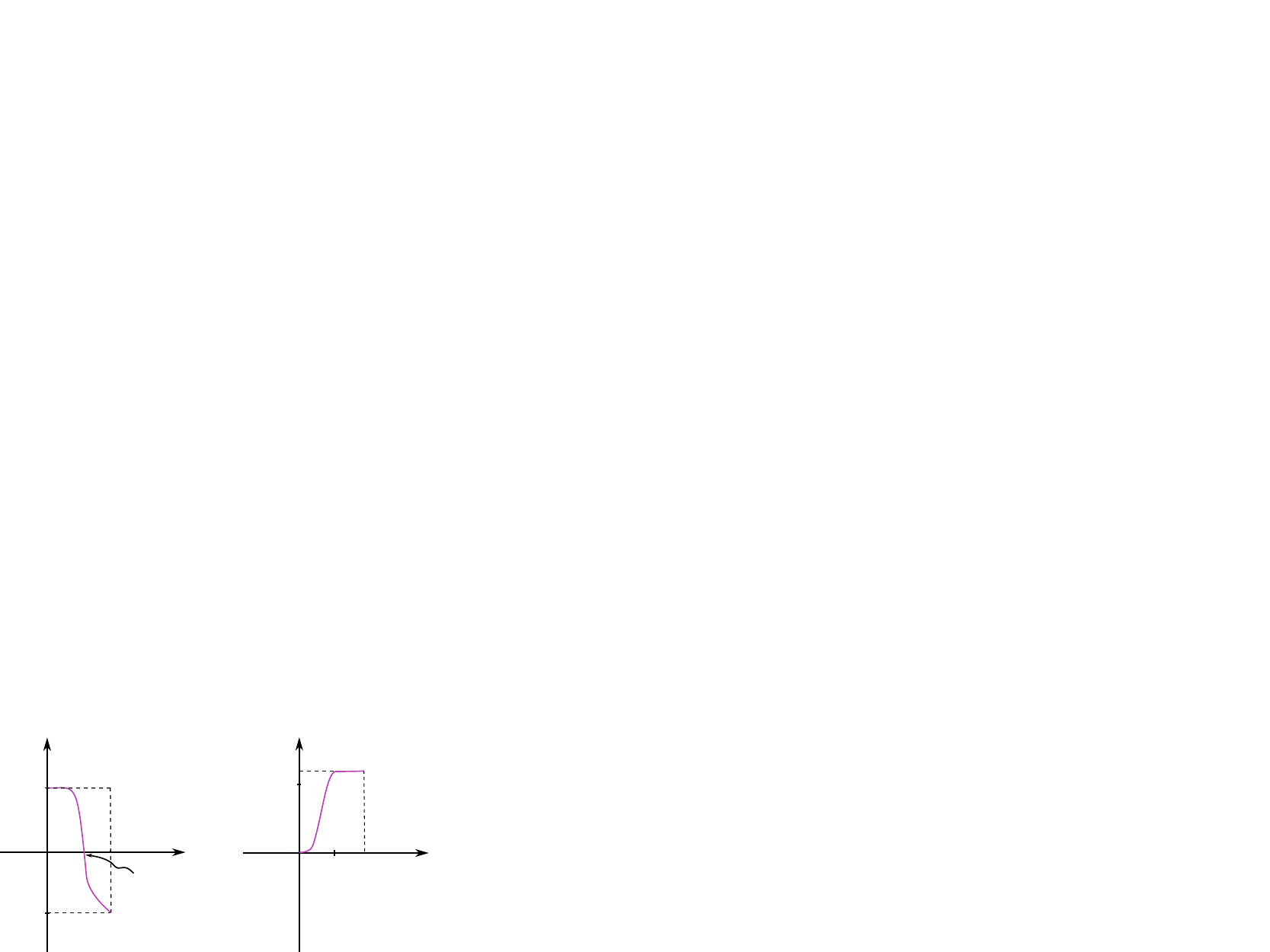}}
\put(-5, 22){$-1-\epsilon$}
\put(55, 80){$h_t$}
\put(203, 80){$g_t$}
\put(265,53){$x$}
\put(110,53){$x$}
\put(170,111){$- \frac{K}{t}$}
\put(22,100){$1$}
\put(89,47){$\delta_{h_t}$}
\put(210,52){$\delta_{g_t}$}
\end{picture}
\caption{$h_t$ and $g_t$ functions}\label{htgt-functions}
\end{center}
\end{figure}

The two families $\{h_t\}_{t\in[-1,0)}$ and $\{g_t\}_{t\in[-1,0)}$ described in Figure
~\ref{htgt-functions} work for our purpose. The conditions ii) and iii) are
satisfied for our choice of $h_t$ and $g_t$, and if $\delta_{h_t}$ and $\delta_{g_t}$ are
such that $h_t(x)<0$ for  $x\in[\delta_{h_t}, 1]$ and $g_t(x)=-\frac{K}{t}$ for
$x\in[\delta_{g_t},1]$, then i) is satisfied as long as
$\delta_{g_t}<\delta_{h_t}$.

Denote the extended family of forms also by $\alpha_t$ and by
$\xi_t$ the family of contact structures given by
$\xi_t=ker(\alpha_t)$, $t\in[-1,0)$. By Gray's theorem there exists a
family of diffeomorphisms $f_t:M\rightarrow M$ such that $(f_t)_*
\xi_{-1} = \xi_t$.

As announced, we would like to construct a family of diffeomorphisms
$\{f_t\}_t$ that fix the pages setwise. Following the proof of
Theorem ~\ref{Gr}, we construct $\{f_t\}_t$ as the flow of a vector
field $X_t\in\xi_t$, for which we have the equality (\ref{eq4}).

We already know that such a $X_t$ exists but would need it to be
tangent to the page. First notice that
$\frac{d\alpha_t}{dt}=\frac{1}{t^2}K d\theta$ and choose some
vector $v\in T\Sigma\cap \xi_t$. Applying both sides of (\ref{eq4}) to
$v$ we get
\begin{equation}
d\alpha_{t}(X_t, v)=\frac{d\alpha_t}{dt}(v_{\alpha_t})\alpha_t(v) -
\frac{d\alpha_t}{dt}(v)
\end{equation}
As $v\in \xi_t =ker(\alpha_t)$ and $v$ has no $\theta$-component,
the last equality is equivalent to
\begin{equation}
d\alpha_t(X_t, v)=0.
\end{equation}

As $d\alpha_t$ is an area form on $\xi_t$, this implies that $X_t$
and $v$ are linearly dependent and therefore $X_t\in T\Sigma\cap
\xi_t$ ($X_t=0$ at singular points, where $T\Sigma$ and $\xi_t$
coincide).

We are now looking at the singularities of $X_t, t\in[-1,0)$. Because for $t_1\ne t_2$, 
$\alpha_{t_1}$ and $\alpha_{t_2}$ differ by a multiple of $d\theta$,
the flowlines given by $X_{t_1}$ and $X_{t_2}$ on a page $\Sigma_{\theta}$ coincide. Therefore,
it suffices to look at the singularities of $X_{t_0}$ for some $t_0\in[-1,0)$.
On
$\Sigma_{\theta}$ there are no negative elliptic singularities away
from the binding as the contact planes and the planes tangent to the
page almost coincide, as oriented plane fields (a negative elliptic
singularity $e$ would require $\xi_e$ and $T_e\Sigma$ to coincide
but have different orientations). Thus, for each $\theta$, all
points on $\Sigma_{\theta}$, except for singularities of $X_{t_0}$ and
stable submanifolds of hyperbolic points, flow in finite time into
an arbitrarily small neighborhood of the binding.
Define $S_{\theta}$ as the set of points on $\Sigma_{\theta}$ that are either singularities of $X_{t_0}$ or on stable submanifolds of hyperbolic points.
Let $\mathcal{S}=\cup S_{\theta}$ as $\theta$ varies from 0 to $2\pi$.

First, note that we can arrange the monodromy map $\phi$ to fix the
singularities on the cutting page, by thinking of $\phi$ as a
composition of Dehn twists away from these points. For isolated
values of $\theta$, $X_{t_0}$ might  exhibit connections between
hyperbolic singularities. With that said, $\mathcal{S}$ has a $CW$ structure
with\\
\hspace*{0.2in}1-skeleton: union of $\{x\}\times [0,1]$ for singular points $x$, and connections between
hyperbolic singularities\\
\hspace*{0.2in}2-skeleton: union of stable submanifolds of hyperbolic singularities.

We want to arrange $K$ in such a way that the bad arcs of $K$ (where
$K$ is not braided) are disjoint from $\mathcal{S}$.  Figure \ref{wrinkleA} depicts a bad arc of $K$. At the point $p$ which lies on the bad arc, $K$ intersects the contact plane $\xi_p$ positively, and it intersects the page negatively (the positive normal vectors to the contact plane and page are shown by arrows). We
introduce wrinkles along $K$ as in Figure \ref{wrinkleA}. See below for the
explicit definition.

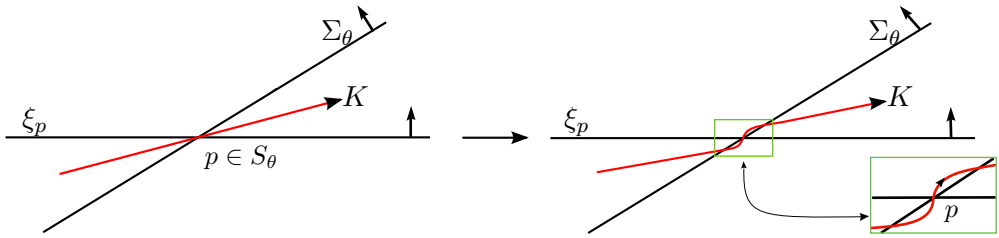
\begin{figure}[htpb]
\centering
\begin{picture}(370, 80)
\put(0,0){\includegraphics{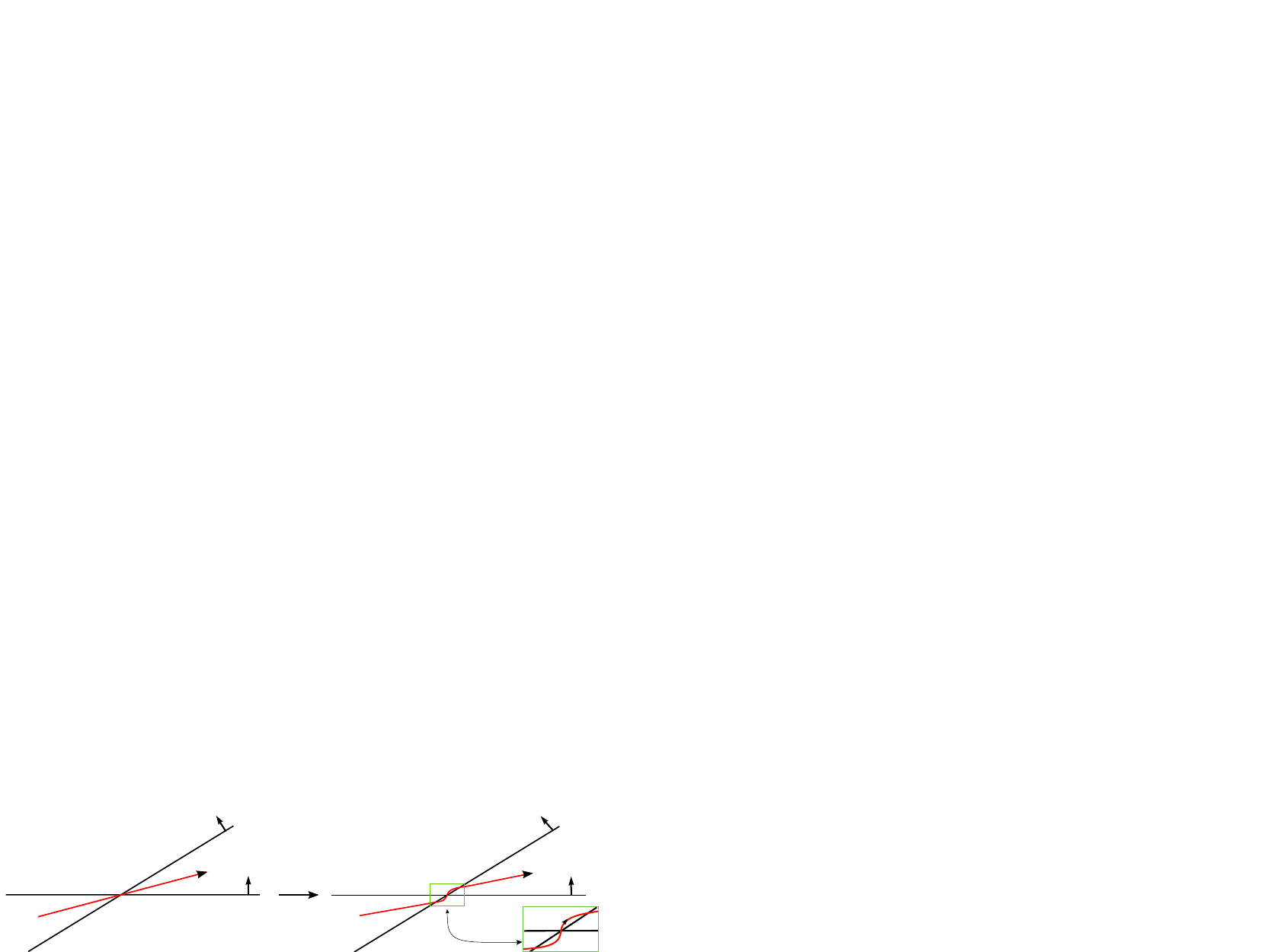}}
\put(131, 48){$K$}
\put(337, 48){$K$}
\put(10, 40){$\xi_p$}
\put(215, 40){$\xi_p$}
\put(79,26){\small{$p\in S_{\theta}$}}
\put(359,6){\small{$p$}}
\put(123, 73){$\Sigma_\theta$}
\put(330, 73){$\Sigma_\theta$}
\end{picture}
\caption{Wrinkling $K$ in order to avoid intersections with
  $\mathcal{S}$}\label{wrinkleA}

\end{figure}

The wrinkles may increase the number of arcs where the link is not
braided but this is fine, as these new arcs avoid $\mathcal{S}$.

By general position, we may assume $K\cap(1-\textrm{skeleton of }
\mathcal{S})=\emptyset$ and $K\pitchfork(2-\textrm{skeleton of } \mathcal{S})$ is a
finite number of points.  A small neighborhood of $p\in
K\pitchfork(2-\textrm{skeleton of } \mathcal{S})$ in $\Sigma_{\theta}$ is
foliated by intervals $(-\epsilon, \epsilon)$, in the same way as a
small disk in the $xy$-plane centered at (0,1,0) in $(\mathbb{R}^3,
\xi_{std})$. It follows that $p$ has a neighborhood in $M$ which
is contactomorphic to a neighborhood of $q=(0,1,0)$ in
$(\mathbb{R}^3, \xi_{std})$. Consider the standard $(x,y,z)$
coordinate system in such a neighborhood of $q$. The contact plane at $q$
is given by the equation $z=x$. As at $p$ the contact
plane and the plane tangent to the page almost coincide, we consider the
plane $z=(1+\epsilon) x$ at $q$ to correspond to the plane tangent to the page at $p$. We identify the bad arc of $K$ with the segment given by $y=1$,
$z=(1+\frac{\epsilon}{2})x$, $x\in [-\delta, \delta]$.  
We call  this segment $W$. See Figure \ref{wrinkleA_coords}. With this setting, the wrinkle takes $W$ to $\tilde{W}$ with the
following properties:\\
\hspace{0.2in}i) $\tilde{W}$ is given by $z=(1+\frac{3\epsilon}{2})x$, $y=1$  for $x\in (-\frac{\delta}{3}, \frac{\delta}{3})$. \\
\hspace{0.2in}ii) $\tilde{W}$ is given by $z=(1+\frac{\epsilon}{2})x$, $y=1$ for $x\in [-\delta, -\frac{2\delta}{3})\cup (\frac{2\delta}{3}, \delta]$
\\
\hspace{0.2in}iii) $\frac{dz}{dx}>0$ along $\tilde{W}$ for $x\in (-\delta, \delta)$.\\\\\\

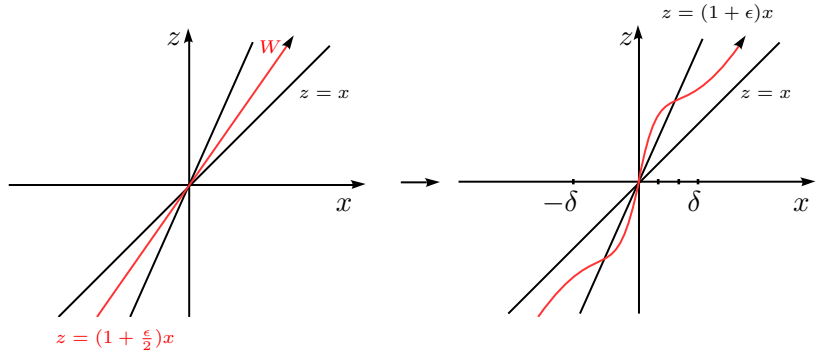
\begin{figure}[htpb]
\centering
\begin{picture}(370, 80)
\put(0,0){\includegraphics{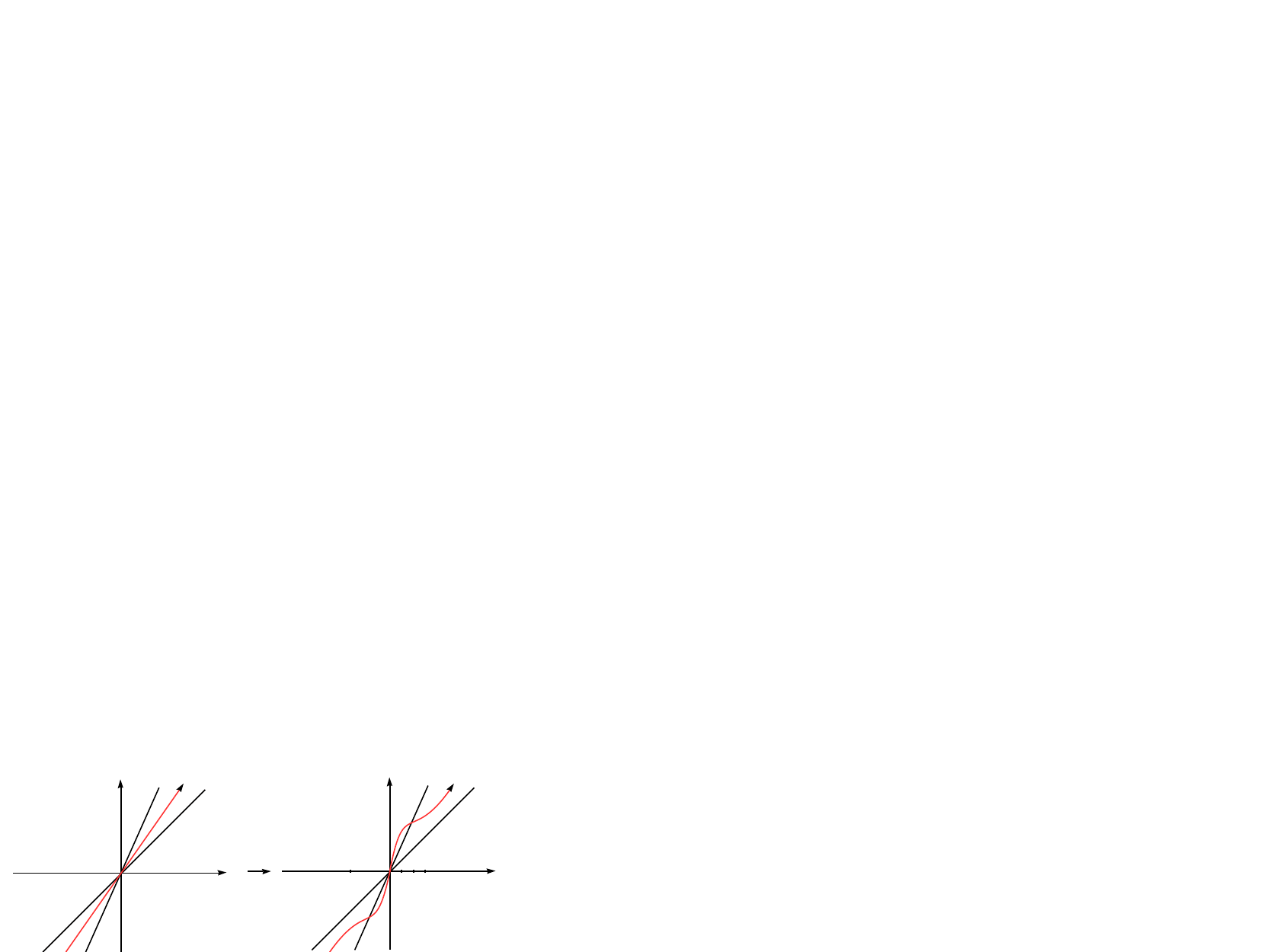}}
\put(132, 40){$\small{x}$}
\put(305, 40){$\small{x}$}
\put(265, 40){$\tiny{\delta}$}
\put(210, 40){$\tiny{-\delta}$}
\put(68, 103){$z$}
\put(239, 103){$z$}
\put(103,100){\color{red}\tiny$W$}
\put(285,83){\tiny$z=x$}
\put(255,113){\tiny$z=(1+\epsilon)x$}
\put(26,-10){\color{red}{\tiny$z=(1+\frac{\epsilon}{2})x$}}
\put(118,83){\tiny$z=x$}
\end{picture}
\caption{Wrinkling $W$.}\label{wrinkleA_coords}

\end{figure}

 Condition iii) guarantees that the link remains transverse throughout the wrinkling.

The arc $\tilde{W}$ still intersects the set $\mathcal{S}$ at the point $q\in \tilde{W}$, but after wrinkling this intersection no longer sits on a bad arc.

We perform this wrinkling at each point where a bad arc of $K$ intersects $\mathcal{S}$, and obtain a link $\tilde{K}$ which does not intersect $\mathcal{S}$ along any of its bad arcs.  We now apply the family of
diffeomorphisms $\{f_t\}_{t\in[-1,0)}$. Then $\tilde{K_{\epsilon}}:= f_{\epsilon}(\tilde{K})$ has
all its  bad arcs in a neighborhood of the binding.  By Theorem
\ref{thm:bennequin}, there is a transverse isotopy taking $\tilde{K_{\epsilon}}$ to a braid $B$. This braid $B$ is transversely isotopic to $K$.

\end{proof}

\section{Markov's Theorem in an Open Book Decomposition}

\begin{figure}[htpb]
\centering
\begin{picture}(390, 80)
\put(0,0){\includegraphics{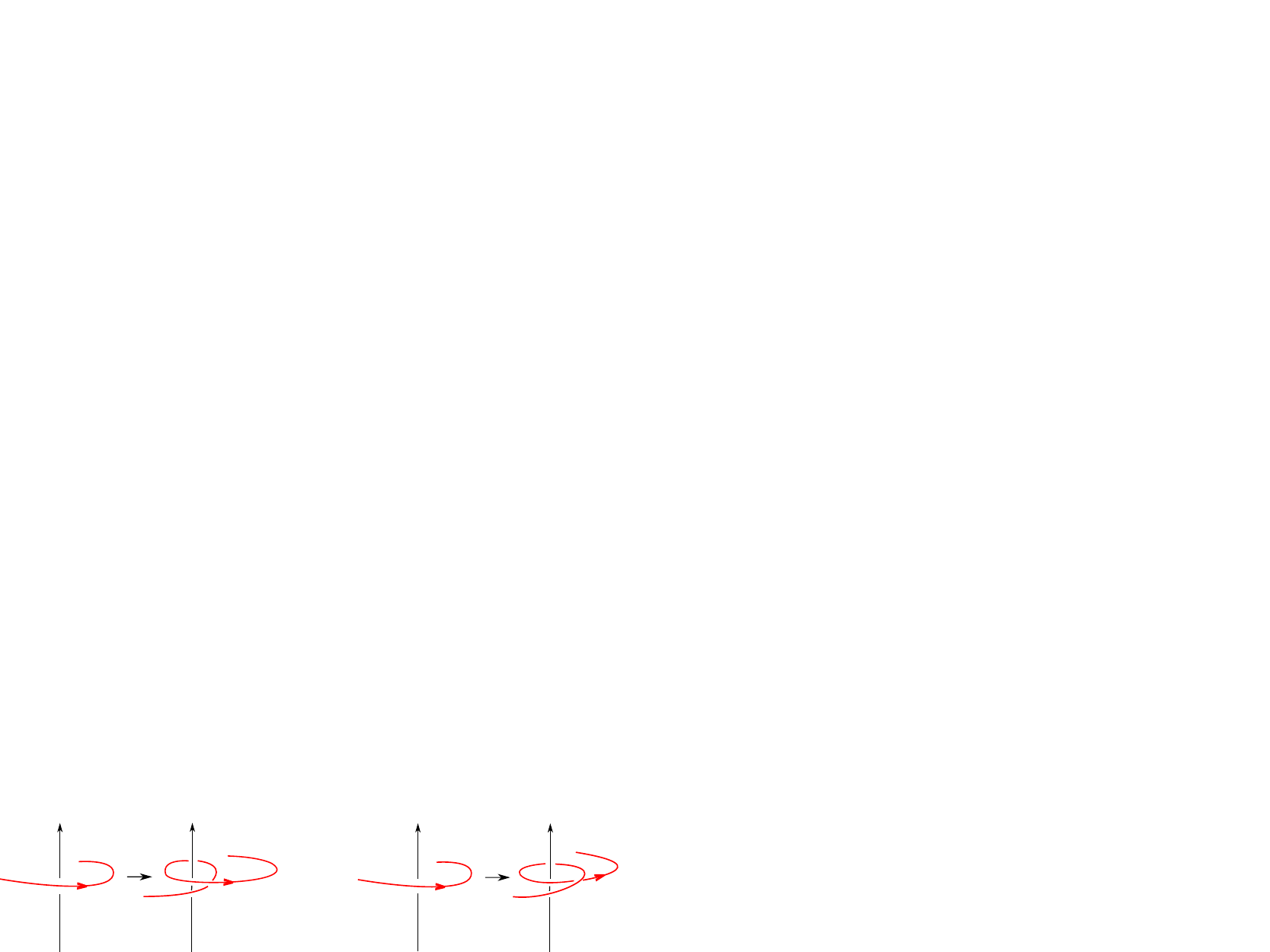}}
\put(27,72){$\small{L}$}
\put(253,72){$\small{L}$}
\end{picture}
\caption{Positive (left) and negative (right) Markov moves.}\label{Markovmove}

\end{figure}

The goal of this section is to generalize the following theorem
\begin{theorem} (Orevkov, Shevchishin \cite{OrSh}) In $(\mathbb{R}^3, \xi_{std})$ two braids
represent transversely isotopic links if and only if one can pass
from one braid to the other by braid isotopies, positive Markov
moves and their inverses. \label{Thm:OrSh}
\end{theorem}
A \emph{Markov move} near a component of the binding $L$ is given by
the introduction of a new loop around $L$, as in Figure \ref{Markovmove}. These
moves are done near the binding, where the standard coordinates
apply. While a positive Markov move (and its inverse) keeps the link
transverse, a negative Markov move (and its inverse) fail to do so
(see \cite{Be} for details).

\begin{theorem} Let $M$ be a 3-dimensional, closed,
oriented manifold and $(L, \pi)$ an open book decomposition for $M$
together with its supported contact structure $\xi$. Let $K_0$ and
$K_1$ be transverse closed braid representatives of the same
topological link. Then $K_0$ and $K_1$ are transversely isotopic if
and only if they differ by braid isotopies and positive Markov moves
and their inverses. \label{Thm:Markovtr}
\end{theorem}
The topological version has been previously proven by Skora
\cite{Sk} and Sundheim \cite{Su}. This case immediately follows from
the proof of the transverse case, as one does not need to worry
about transversality throughout the isotopy. If transversality is
not required, both positive and negative Markov moves are
allowed.
\begin{proof}
The reverse implication is straightforward. One needs to note
that an isotopy through braids is done away from the binding. As the
contact planes almost coincide with the planes tangent to the pages
this isotopy stays transverse with respect to the contact structure.

The direct implication takes more work. Let $K_0$ and $K_1$ be
transverse braid representatives of the same topological link $K$
and $\{K_t\}_{t\in[0,1]}$ a transverse isotopy from $K_0$ to $K_1$.
We parametrize the isotopy by $\mathcal{K}:\coprod S^1\times
[0,1]\rightarrow M$, such that $\mathcal{K}_t$ defined by
$s\rightarrow \mathcal{K}(s,t)$ is a parametrization of $K_t$, where
$s$ is the positively oriented coordinate on each $S^1$. Denote by $\mathfrak{A}\subset M$ the immersed annulus $\mathcal{K}(\coprod S^1\times [0,1])$.

Let $\theta$ be the positive coordinate normal to the page. A
\textit{bad zone} of $\mathcal{K}$ is a connected component of  $\mathfrak{A}$ where
 $\frac{\partial}{\partial{s}}(\pi\circ\mathcal{K})\le0$.
 We denote by $\mathcal{B}$ the union of all bad zones of $\mathcal{K}$.

 We would like to take all the bad zones of $\mathcal{K}$ in a
neighborhood of the binding. This way, the proof is reduced to the
standard case proved by Orevkov and Shevchishin in \cite{OrSh}. For
this, we use the family of diffeomorphisms $\{f_t\}_{t\in [-1,0)}$
constructed in the proof of Theorem \ref{Alexander}. 
Keeping the same notations, we want to arrange the isotopy $\mathcal{K}$ in such a way that $\mathcal{B}\cap \mathcal{S}=\emptyset$.

By general position, $\mathcal{B}\cap \mathcal{S}$ consists of arcs $l$, along which  $\frac{\partial(\pi\circ\mathcal{K})}{\partial s}\le0$ and points at which $\frac{\partial(\pi\circ\mathcal{K})}{\partial s}\le0$. 
The arcs $l$ are intersections of the 2-cells of $\mathcal{S}$ with $\mathcal{B}$ while the points are intersections of the 1-cells of $\mathcal{S}$ with $\mathcal {B}$.  
We modify the isotopy $\mathcal{K}$ in such a way so as to replace the arc $l\subset \mathcal{B}\cap \mathcal{S}$ by another arc $l'\subset \mathcal{S}$, along which $\frac{\partial(\pi\circ\mathcal{K})}{\partial s}>0$, and to eliminate the isolated points of $\mathcal{B}\cap\mathcal{S}$.  We describe the process below.

Let  $l\subset \mathcal{B}\cap \mathcal{S}$ and identify a small region of $\mathcal{B}$ containing $l$ with the rectangle $R = [\frac{1}{2}, \frac{3}{2}] \times [0,1]$ (denote by $t$ the coordinate on the first factor) such that $l$ is identified with $[\frac{1}{2}, \frac{3}{2}]\times \{\frac{1}{2}\}$ and let $R_t := \{t\}\times [0,1]\subset R$, $t \in [\frac{1}{2}, \frac{3}{2}]$.  We make the following identifications: 
\begin{itemize}
\item For  $t \in [\frac{1}{2}, \frac{3}{2}]$ identify $R_t$ with the segment $W_t\subset \mathbb{R}^3$ given by $y=t$ and $z= ( t+\frac{\epsilon}{2})x$.
\item For  $t \in [\frac{1}{2}, \frac{3}{2}]$ identify the contact plane at $P(t, \frac{1}{2})\in l$ with the contact plane at $(0,t,0)$, that is the plane given by $z= tx$.
\item For  $t \in [\frac{1}{2}, \frac{3}{2}]$ identify the page of the open book at  $(t, \frac{1}{2})\in l$ with the plane at $(0,t,0)$ given by $z= (t+\epsilon)x$.\\\\
\end{itemize}

\begin{figure}[htpb]
\centering
\begin{picture}(370, 80)
\put(0,0){\includegraphics{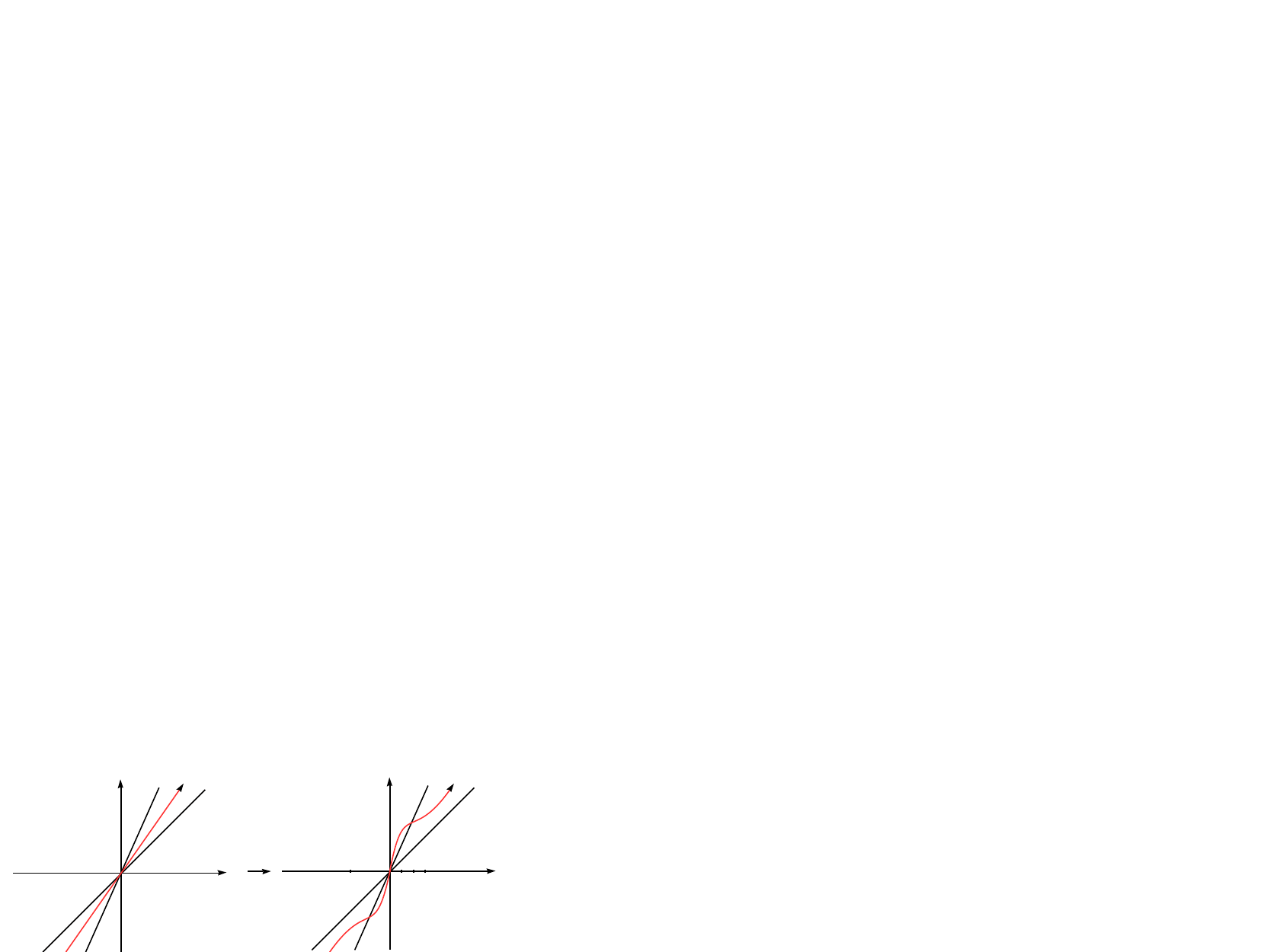}}
\put(132, 40){$\small{x}$}
\put(305, 40){$\small{x}$}
\put(265, 40){$\tiny{\delta}$}
\put(210, 40){$\tiny{-\delta}$}
\put(68, 103){$z$}
\put(239, 103){$z$}
\put(115,100){\color{red}\tiny$W_t$}
\put(205,-10){\color{red}{\tiny$\tilde{W_t}$}}
\put(285,83){\tiny$z=tx$}
\put(255,113){\tiny$z=(t+\epsilon)x$}
\put(26,-10){\color{red}{\tiny$z=(t+\frac{\epsilon}{2})x$}}
\put(118,83){\tiny$z=tx$}
\end{picture}
\caption{Wrinkling $W_t$.}\label{wrinkleM_coords}

\end{figure}

\begin{figure}[htpb]
\begin{center}
\begin{picture}(396, 144)
\put(0,0){\includegraphics{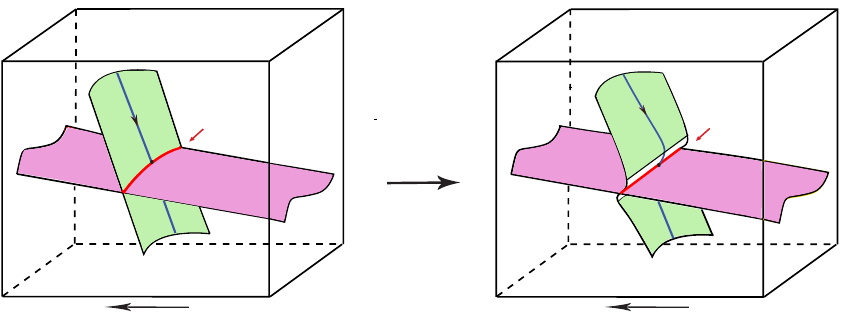}}
\put(-5, 7){\tiny{$\Sigma\times \theta_1$}}
\put(115, 7){\tiny{$\Sigma\times \theta_0$}}
\put(230, 7){\tiny{$\Sigma\times \theta_1$}}
\put(355, 7){\tiny{$\Sigma\times \theta_0$}}
\put(70,5){\small{$\theta$}}
\put(310,5){\small{$\theta$}}
\put(99, 100){\color{red}$l$}
\put(343, 100){\color{red}$l'$}
\put(67,47){$\mathcal{B}$}
\put(310, 47){$\mathcal{B}$}
\put(12,80){\small{$\mathcal{S}$}}
\put(250,80){\small{$\mathcal{S}$}}
\put(75,80){\small{$x$}}
\put(320,80){\small{$x$}}
\end{picture}
\caption{Isotopy modification along an $l$ arc. The new intersection arc $l'$ sits outside of $\mathcal{B}$  ( green region) }\label{im}
\end{center}
\end{figure}

With this settings, the wrinkle takes $W_t$ to $\tilde{W_t}$, $t \in [\frac{1}{2}, \frac{3}{2}]$, with the
following properties:\\
\hspace{0.2in}i) $\tilde{W_t}$ is given by $z=(t+\frac{3\epsilon}{2})x$, $y=t$  for $x\in (-\frac{\delta}{3}, \frac{\delta}{3})$. \\
\hspace{0.2in}ii) $\tilde{W}$ is given by $z=(t+\frac{\epsilon}{2})x$, $y=t$ for $x\in [-\delta, -\frac{2\delta}{3})\cup (\frac{2\delta}{3}, \delta]$
\\
\hspace{0.2in}iii) $\frac{dz}{dx}>0$ along $\tilde{W_t}$ for $x\in (-\delta, \delta)$.
See Figure \ref{wrinkleM_coords}.
\\

We perform this wrinkling along every arc in the intersection $ \mathcal{B}\cap \mathcal{S}$. If two such arcs $l_1$ and $l_2$ intersect each other there exist two disjoint regions on $\mathfrak{A}(l_1\cup\l_2)$ one containing $l_1$ and the other containing $l_2$ which we can identify with rectangles $R_1$ and $R_2$ as above. We wrinkle along $l_1$ and $l_2$ separately.

The wrinkling translates to the original setup as in Figure \ref{im}. The bad zone (colored green in Figure \ref{im}) no longer intersects $\mathcal{S}$.

For an isolated point $x\in\mathcal{B}\cap \mathcal{S}$, we pick a small bad arc $l$ containing $x$ on the isotopy annulus and we do the wrinkling along $l$ as above. 

We obtain an isotopy annulus which does not intersect $\mathcal{S}$ along any bad arcs. 
\end{proof}

\end{document}